\documentclass[technote]{IEEEtran}
\ifCLASSINFOpdf
\else
\fi
\usepackage{graphicx} 
\usepackage{epsfig} 
\usepackage{amsmath} 
\usepackage{amssymb}  
\usepackage{psfrag}
\usepackage{color}
\usepackage{cite}

\graphicspath{{Figures/}}

\newcommand{\tA}{\mathcal{A}}
\newcommand{\tB}{\mathcal{B}}

\newcommand{\tP}{\mathcal{P}}

\newcommand{\tF}{\mathcal{F}}
\newcommand{\di}{\mathrm{diag}}
\newcommand{\tG}{\mathcal{G}}

\newtheorem{thm}{Theorem} 
 
\newtheorem{lem}{Lemma}

\newtheorem{defi}{Definition}
\newtheorem{rem}{Remark}

\begin{document}
%
\title{Robust Consensus Analysis and Design under Relative State Constraints or Uncertainties}
%
%
%
\author{Dinh~Hoa~Nguyen,~\IEEEmembership{Member,~IEEE,}
        Tatsuo~Narikiyo, 
        and~Michihiro~Kawanishi,~\IEEEmembership{Member,~IEEE}
\thanks{The authors are with Control System Laboratory, Department of Advanced Science and Technology, Toyota Technological Institute, 
        2-12-1 Hisakata, Tempaku-ku, Nagoya 468-8511, Japan. Emails: 
        N.DinhHoa.Nguyen@ieee.org, n-tatsuo@toyota-ti.ac.jp, kawa@toyota-ti.ac.jp}}

\newcounter{MYtempeqncnt}

\maketitle


\begin{abstract}
 
This paper proposes a new approach to analyze and design distributed robust consensus control protocols for general linear leaderless multi-agent systems (MASs) 
in presence of relative-state constraints or uncertainties.  
First, we show that the MAS robust consensus under relative-state constraints or uncertainties is equivalent to the robust stability under state constraints or uncertainties of a transformed MAS. 
Next, the transformed MAS under state constraints or uncertainties is reformulated as a network of Lur'e systems. 
By employing S-procedure, Lyapunov theory, and Lasalle's invariance principle, a sufficient condition for robust consensus and the design of robust consensus controller gain are derived from solutions of a distributed LMI convex problem. 
Finally, numerical examples are introduced to illustrate the effectiveness of the proposed theoretical approach.

\end{abstract}



\section{Introduction}

Multi-agent systems (MASs) and their cooperative control problems have been extensively studied and applied to many practical systems, e.g.,     
power grids, wireless sensor networks, transportation networks, systems biology, etc., because of their key advantage of achieving 
global objectives by performing local measurements and controls at each agent and simultaneously collaborating among agents using that local information. 
Among many interesting problems, consensus is one of the most important and intensively investigated issues in MASs 
due to its attraction in both theory and applications \cite{Olfati-Saber:2004,Olfati-Saber:2007,Ren:2007}.

In practical MASs, agents' inputs or states and the exchanged information among agents are subjected to constraints or uncertainties 
due to physical limitations of agents or uncertain communication channels. 
Realistic examples are consensus of vehicles with limited speeds and working space, smart buildings energy control with temperature and humidity are required in specific ranges, 
just to name a few. 
Therefore, the MAS consensus under constraints and uncertainties on the inputs, states, or relative states of agents is a significant, realistic problem and is worth studying.  
However, this problem was not investigated in the early researches on MASs and 
it just has been considered in some recent studies \cite{Nedic:2010,Lin:2014,Lee:2011,Liu:2012,Wang:2013,Meng:2012,Su:2013,Takaba:2014,Zhang:2014,Zhang:2014j,Lim:2015}.  

A constrained consensus problem was investigated in \cite{Nedic:2010} where the states of agents are required to lie in individual closed convex sets and the final consensus state must belong to the non-empty intersection of those sets. Accordingly, a projected consensus algorithm was proposed and then applied to distributed optimization problems. Following this research line, \cite{Lin:2014} extended the result in \cite{Nedic:2010} to the context where communication delays exist.  
In another work, \cite{Lee:2011} studied the state increment by utilizing the model predictive control (MPC) method. 
However, distributed and fast MPC algorithms need to be further developed in order to use in large-scale MASs. 
Another direction is to employ the discarded consensus algorithms \cite{Liu:2012,Wang:2013}. 
Nevertheless, a requirement of these approaches as well as in \cite{Nedic:2010,Lin:2014} is that the initial states of agents 
must belong to some sets specified by the constraints, i.e., the consensus is only local. 
Moreover, only agents with single integrator dynamics were considered in \cite{Liu:2012,Wang:2013}.  

To achieve the global or semi-global consensus in presence of input or state constraints, some consensus laws were presented in \cite{Meng:2012,Su:2013}, but they were only for leader-follower MASs. 
In other researches, \cite{Takaba:2014,Zhang:2014,Zhang:2014j} derived global consensus under input or state constraints by reformulating the constrained MAS as a network of Lur'e systems and employing Lyapunov theory. 
The paper \cite{Takaba:2014} considered linear agents with input saturation but agents' dynamics is limited to be single-input. 
Next, \cite{Zhang:2014} and \cite{Zhang:2014j} investigated consensus problems for general linear MASs where outputs of agents are incrementally bounded or passive 
and obtained sufficient conditions for global consensus in the form of LMI convex problems.   

On the other hand, the MAS consensus under relative-state constraints has been recently studied in \cite{Lim:2015} 
within a very special context where the input matrices of agents are identity matrices and the consensus controller gain is a diagonal matrix. 
Then sufficient conditions were proposed for the cases of $2$-norm and $\infty$-norm bounded constraints on relative states of agents. 
Nevertheless, the consensus is only local and no consensus controller design was given in \cite{Lim:2015}. 

This paper proposes a new approach to analyze and design {\it distributed robust consensus controllers} for {\it general linear homogeneous leaderless MASs} to achieve {\it global consensus} 
under {\it relative-state constraints or uncertainties} which are in the form of a {\it sector-bounded condition}. 
Our approach covers broader systems and scenarios than those in the existing researches, and hence constitutes our first contribution. 
Consequently, we further develop the edge dynamics proposed in \cite{Nguyen:2015TAC1} to achieve that the currently considered problem is {\it equivalent to a distributed robust stabilization problem under state constraints or uncertainties for a transformed MAS}.  This serves as our second contribution.  
Next, the transformed MAS is rewritten as a network of Lur'e systems and the robust stabilization problem is formulated as a distributed convex LMI problem.   
In comparison with the one in \cite{Takaba:2014} for a similar type of Lur'e networks, our LMI problem is less conservative that:  
(i) employs a more general method namely the S-procedure; 
(ii) gives an exponential convergence to consensus instead of asymptotic convergence.  
Furthermore, our consensus controller gain is much more general than the diagonal one in \cite{Lim:2015}. 
Those advantages clearly show our third contribution. 




The following notation and symbols will be used in the paper. $\mathbb{R}$ and $\mathbb{C}$ stand for the real 
and complex sets. 
Moreover, $\mathbf{1}_n$ denotes the $n\times1$ vector 
with all elements equal to $1$, and $I_n$ denotes the $n\times n$ identity matrix.
Next, $\otimes$ stands for the Kronecker product, $\di\{\}$ denotes diagonal or block-diagonal matrices, and   
$\mathrm{sym}(A)$ denotes $A+A^T$ for any real matrix $A$.
Lastly, $\succ$ and $\succeq$ denote the positive definiteness and positive semi-definiteness of a matrix, and similar meanings are used for $\prec$ and $\preceq$.

\section{Problem Description}
\label{prob}

Consider a MAS consisting of $N$ identical agents with the following linear dynamics
\begin{equation}
	\label{agent}
		\dot{x}_i = Ax_i+Bu_i,  i=1,\ldots,N,
\end{equation} 
where $x_i \in \mathbb{R}^{n}$ is the state vector, $u_i \in \mathbb{R}^{m}$ is the control input,  
$A \in \mathbb{R}^{n\times n}$, $B \in \mathbb{R}^{n \times m}$. 
The whole MAS is then described by
\begin{equation}
	\label{mas}
		\dot{x} = (I_{N}\otimes A)x+(I_{N}\otimes B)u,
\end{equation}
where $x = \left[x_{1}^T,\ldots,x_{N}^T\right]^T, u=\left[u_{1}^T,\ldots,u_{N}^T\right]^T$. 
Let $\tG$ be an undirected graph representing the information structure in the MAS, 
in which each node in $\mathcal{G}$ represents an agent and each edge in $\mathcal{G}$ represents 
the interconnection between two agents. 
Denote $L \in \mathbb{R}^{N \times N}$ and $E \in \mathbb{R}^{N\times M}$ the Laplacian matrix and the incidence matrix associated with $\tG$. 
Then $L=EE^T$ and $E^T\mathbf{1}_{N}=0$. 

The following assumptions will be employed. 
\begin{itemize}
	\item[\bf A1:]  $(A,B)$ is stabilizable. 
	\item[\bf A2:]  All eigenvalues of $A$ is on the closed left half complex plane. 
	\item[\bf A3:]  $\tG$ is undirected and connected.  
\end{itemize}
Assumptions A1--A2 are necessary and sufficient such that the consensus can be achieved and stable (see e.g. \cite{Ma:2010}). 
Next, the consensus of agents is defined as follows. 
\begin{defi}
The MAS with linear dynamics of agents represented by (\ref{agent}) and the information exchange among agents represented by 
$\mathcal{G}$ is said to reach a consensus if 
\begin{equation}
	\label{consensus-defi}
	\lim_{t \rightarrow \infty} \|x_{i}(t)-x_{j}(t)\| = 0 ~\forall~i,j=1,\ldots,N.
\end{equation}
\end{defi}
Due to physical limitations on the communication range and bandwidth of agents or uncertain information channels, 
the exchanged relative states among agents could be bounded or contain some uncertainties. 
To take into account those practical issues in the control analysis and design, 
we define in the following a new state vector and a new control input, $$z\triangleq (E^T\otimes I_{n})x, ~w\triangleq (E^T\otimes I_{m})u.$$ 
Let $L^{\dag}$ be the generalized pseudoinverse of the Laplacian matrix $L$ \cite{Gutman:2004}. 
Then multiplying both sides of (\ref{mas}) with $E^T\otimes I_{n}$ gives us the following edge dynamics \cite{Nguyen:2015TAC1},
\begin{equation}
	\label{edge-dyn}
	\dot{z} = [(E^TL^{\dag}E)\otimes A]z+(I_{M}\otimes B)w.
\end{equation}
Since $z$ composes of all relative states of connected agents, (\ref{mas}) is consensus if this edge dynamics is stabilized. 
Moreover, all relative-state constraints and uncertainties are now represented in term of $z$. 
Accordingly, the following control scenario shall be investigated. 
\begin{itemize}
	\item {\bf Relative-State Constraints/Uncertainties:} For all $j \in [1,M]$, $y_{j,k}=\phi_{j,k}(z_{j,k}) \,\forall\, k=1,\ldots,n$ where $y_{j,k} \in \mathbb{R}$ is the $k$th component of the 
	signal exchanged through the edge $j$; $\phi_{j,k}: \mathbb{R} \rightarrow \mathbb{R}$ is a continuous function that satisfies the following sector-bounded condition:
		\begin{equation}
			\label{rel-state-cstrt}
			\begin{aligned}
				&(\phi_{j,k}(z_{j,k})-\sigma_{k,1}z_{j,k})(\phi_{j,k}(z_{j,k})-\sigma_{k,2}z_{j,k}) \leq 0 \\ 
				& \, \forall \, k=1,\ldots,n; \, \forall \, j=1,\ldots,M,
			\end{aligned}
		\end{equation}
		where $\sigma_{k,1},\sigma_{k,2} \in \mathbb{R}$ are known constants, $\sigma_{k,1} < \sigma_{k,2}$. 		
\end{itemize}

Consequently, we present the control analysis and design problem considered in this paper. 
\begin{itemize}
	\item {\bf Global robust consensus under relative-state constraints or uncertainties:} For the given linear MAS with dynamics of agents represented by (\ref{agent}) and the information 
	exchange among agents represented by $\tG$, find a condition and a control strategy to achieve consensus of agents in the sense of (\ref{consensus-defi}) subjected to the relative-state constraints or 
	uncertainties (\ref{rel-state-cstrt}), for any initial conditions of agents.  
\end{itemize}

\section{Consensus Analysis and Design under Relative-State Constraints or Uncertainties}
\label{rel-state}

\subsection{Equivalence to Robust Stabilization in Presence of State Constraints or Uncertainties}

Denote $\bar{L}=E^TL^{\dag}E$ and $L_{e}=E^TE$.

\begin{lem} \cite{Nguyen:2015TAC1}
\label{L-properties}
The following statements hold.
\begin{itemize}
	\item[(i)] $L_{e}$ has exactly $N-1$ non-zero eigenvalues, which are equal to positive eigenvalues of $L$ while all other eigenvalues of $L_{e}$ if exist are $0$. 
	\item[(ii)] $\bar{L}$ has exactly $N-1$ non-zero eigenvalues, which are all equal to $1$, and other eigenvalues of $\bar{L}$ if exists are $0$.
\end{itemize}
\end{lem} 

Let $U\in\mathbb{R}^{M\times M}$ be an orthogonal matrix that diagonalizes $\bar{L}$, and 
$$\tilde{z} \triangleq (U^T\otimes I_{n})z,\tilde{w} \triangleq (U^T\otimes I_{m})w.$$ 
Subsequently, we obtain from (\ref{edge-dyn}) that
\begin{equation}
	\label{edge-dyn-1}
	\dot{\tilde{z}} = \left(\bar{\Gamma}\otimes A\right)\tilde{z}+(I_{M}\otimes B)\tilde{w},
\end{equation} 
where $\bar{\Gamma}=\di\{0, I_{N-1}\}$ includes all eigenvalues of $\bar{L}$ in its diagonal (due to Lemma \ref{L-properties}). 
Now, let us partition $U$, the state and input vectors in (\ref{edge-dyn-1}) as follows,
\begin{equation}
	U = \begin{bmatrix} U_1 & U_2 \end{bmatrix},
	\tilde{z} = \begin{bmatrix} \tilde{z}_1 \\ \tilde{z}_2 \end{bmatrix},
	\tilde{w} = \begin{bmatrix} \tilde{w}_1 \\ \tilde{w}_2 \end{bmatrix},
\end{equation}
where $U_1\in\mathbb{R}^{M\times(M-N+1)}$, $U_2\in\mathbb{R}^{M\times(N-1)}$, $\tilde{z}_1\in\mathbb{R}^{n(M-N+1)}$, $\tilde{z}_2\in\mathbb{R}^{n(N-1)}$, $\tilde{w}_1\in\mathbb{R}^{n(M-N+1)}$, 
$\tilde{w}_2\in\mathbb{R}^{n(N-1)}$. Then (\ref{edge-dyn-1}) is equivalent to
\begin{equation}
	\label{edge-dyn-2}
	\begin{aligned}
		\dot{\tilde{z}}_1 &= (I_{M-N+1}\otimes B)\tilde{w}_1, \\
		\dot{\tilde{z}}_2 &= (I_{N-1}\otimes A)\tilde{z}_2+(I_{N-1}\otimes B)\tilde{w}_2, 
	\end{aligned}
\end{equation}
and $\tilde{z}_1=(U_1^T\otimes I_{n})z$, $\tilde{z}_2=(U_2^T\otimes I_{n})z$. 

Let $\Gamma \in \mathbb{R}^{(N-1)\times(N-1)}$ be the diagonal matrix including all non-zero eigenvalues of $L$ in its diagonal, and $V \in \mathbb{R}^{N\times N}$ is an orthogonal matrix such that 
\begin{equation}
	V^TLV = \begin{bmatrix} 0 & 0 \\ 0 & \Gamma \end{bmatrix}.
\end{equation}
Partitioning $V$ into $[V_1,V_2]$ where $V_1 \in \mathbb{R}^{N}$, $V_2 \in \mathbb{R}^{N\times(N-1)}$. Then 
\begin{equation}
	\label{V2-eq}
	LV_2=V_2\Gamma \Leftrightarrow V_2^TLV_2=\Gamma,
\end{equation}
since $V_2^TV_2=I_{N-1}$. 

Denote
$\Phi(z) \triangleq [\phi_{1}^T(z_{1}),\ldots,\phi_{M}^T(z_{M})]^T$, $\phi_{j}(z_{j}) \triangleq [\phi_{j,1}^T(z_{j,1})$, $\ldots,\phi_{j,N}^T(z_{j,M})]^T$, $\forall\; z=[z_{1}^T,\ldots,z_{M}^T]^T$.  

\begin{thm}
\label{ctlr-equi-cstrt}
Let $U_2$ be chosen as $E^TV_2\Gamma^{-1/2}$, then the following distributed robust stabilizing controller in presence of state constraints or uncertainties
\begin{equation}
	\label{stab-ctlr-cstrt}
	\tilde{w}=\tF \left(U^{T} \otimes I_{n}\right) \Phi(z),
\end{equation}
for the transformed edge dynamics (\ref{edge-dyn-2}) with 
$\tF=F \otimes K$, $K \in \mathbb{R}^{m\times n}$ and 
\begin{equation}
	\label{F-form}
	F = \begin{bmatrix} 0 & 0 \\ 0 & \Gamma \end{bmatrix},
\end{equation}
is equivalent to the following distributed robust consensus controller under relative-state constraints or uncertainties (\ref{rel-state-cstrt}),
\begin{equation}
	\label{rel-state-consensus-ctlr-cstrt}
	u=(E \otimes K)\Phi(z), 
\end{equation}
for the initial MAS (\ref{mas}). Furthermore, $\tilde{z}_1(t)=0 ~\forall~ t\geq 0.$
\end{thm}

\begin{IEEEproof}
First, we show that the orthogonality of $U$ is satisfied with $U_2$ chosen to be $E^TV_2\Gamma^{-1/2}$. 
Indeed, $U_2^TU_1=\Gamma^{-1/2}V_2^TEU_1=0$ since $EU_1=0$ due to a fact that $E\bar{L}=E$. 
Moreover, $U_2^TU_2=\Gamma^{-1/2}V_2^TEE^TV_2\Gamma^{-1/2}=\Gamma^{-1/2}\Gamma\Gamma^{-1/2}=I_{N-1}$. 
 
Next, multiplying to the left of (\ref{V2-eq}) with $E^T$ gives us
\begin{align}
	\label{eq:}
	L_{e}E^TV_2 &= E^TV_2\Gamma, \nonumber \\
	\Leftrightarrow L_{e}E^TV_2\Gamma^{-1/2} &= E^TV_2\Gamma^{1/2}, \nonumber \\
	\Leftrightarrow \Gamma^{-1/2}V_2^TEL_{e}E^TV_2\Gamma^{-1/2} &= \Gamma^{-1/2}V_2^TEE^TV_2\Gamma^{1/2}, \nonumber \\
	\Leftrightarrow U_2^TL_{e}U_2 &= \Gamma. 
\end{align}
On the other hand, $\bar{L}U_1=0$, which leads to $L_{e}U_1=0$ since $L_{e}\bar{L}=L_{e}$. Therefore, we obtain
\begin{equation}
	\label{F}
	F=U^TL_{e}U ~ \Leftrightarrow ~ UFU^T = L_{e}.
\end{equation}
Consequently, 
\begin{align}
	\label{eqq}
	\tilde{w} &= \tF \left(U^{T} \otimes I_{n}\right) \Phi(z), \nonumber \\
	\Leftrightarrow  [(U^TE^T)\otimes I_{m}]u &= \left[\left(FU^T\right) \otimes K\right] \Phi(z).
\end{align}
Since $U^TE^T = [EU_1, EU_2]^T = [0,V_2\Gamma^{1/2}]^T$ and $FU^T=[0,E^TV_2\Gamma^{1/2}]^T$, (\ref{eqq}) is equivalent to
\begin{align}
	\label{eqq-1}
	\left[\left(\Gamma^{1/2}V_2^T\right)\otimes I_{m}\right]u &= \left[\left(\Gamma^{1/2}V_2^TE\right) \otimes K\right] \Phi(z), \nonumber \\
	\Leftrightarrow [(V_{2}V_{2}^T)\otimes I_{m}] u &= [(V_{2}V_{2}^TE) \otimes K]\Phi(z),
\end{align}
by multiplying both to the left and to the right of (\ref{eqq-1}) with $(V_2\Gamma^{-1/2})\otimes I_{m}$. 
Note that $V_{1}=\frac{1}{\sqrt{N}}\mathbf{1}_{N}$, then $V_{2}V_{2}^T=I_{N}-V_{1}V_{1}^T=I_{N}-\frac{1}{N}\mathbf{1}_{N}\mathbf{1}_{N}^T$. Hence, we obtain 
$[(I_{N}-\frac{1}{N}\mathbf{1}_{N}\mathbf{1}_{N}^T)\otimes I_{m}]u = [([I_{N}-\frac{1}{N}\mathbf{1}_{N}\mathbf{1}_{N}^T]E)\otimes K]\Phi(z)=(E \otimes K)\Phi(z)$, since $\mathbf{1}_{N}^TE=0.$ 
This is equivalent to $u=(E \otimes K)\Phi(z)+(\mathbf{1}_{N}\otimes I_{m})u_{0}$ for any $u_{0} \in \mathbb{R}^{m}$. Since we are not interested in self-feedback inputs for agents, $u_{0}=0$ 
or equivalently $$u=(E \otimes K)\Phi(z).$$ 
On the other hand, we have $$\tilde{z}_1=U_1^Tz=[(U_1^TE^T)\otimes I_n]x=0  ~\forall~ t\geq 0,$$ since $EU_1=0$. 
\end{IEEEproof}

Employing the result of Theorem \ref{ctlr-equi-cstrt} to the transformed edge dynamics (\ref{edge-dyn-2}), it can be deduced that we only need to design a distributed robust stabilizing controller for the subsystem
\begin{equation}
	\label{edge-dyn-3}
	\dot{\tilde{z}}_2 = (I_{N-1}\otimes A)\tilde{z}_2+(I_{N-1}\otimes B)\tilde{w}_2, 
\end{equation}
having the following form
\begin{equation}
	\label{stab-ctlr-cstrt-1}
	\tilde{w}_2 = \left[\left(\Gamma U_{2}^T\right) \otimes K\right]\Phi(z),
\end{equation}
which is directly calculated from (\ref{stab-ctlr-cstrt}). 
Thus, the interesting result of Theorem \ref{ctlr-equi-cstrt} is that the distributed robust consensus design (\ref{rel-state-consensus-ctlr-cstrt}) under relative-state constraints or uncertainties for the initial MAS (\ref{mas}) is equivalent to a simpler problem of synthesizing a distributed robust stabilizing controller (\ref{stab-ctlr-cstrt-1}) under state constraints or uncertainties for a new MAS 
(\ref{edge-dyn-3}) which has lower dimension. 
In the next section, we will present an approach to design such a distributed robust stabilizing controller.

\begin{rem}
If $\mathcal{G}$ is a spanning tree then $M=N-1$ and hence $\bar{L}=I_{N-1}$. Then we do not need the additional transformation (\ref{edge-dyn-1}).  
Therefore, all results here and in subsequent sections are derived with $\tilde{w}_{2}$ and $\tilde{z}_{2}$ replaced by $w$ and $z$, respectively. 
\end{rem}

\subsection{Distributed Robust Stabilizing Controller Synthesis}

The transformed edge dynamics (\ref{edge-dyn-3}) together with the robust stabilizing controller (\ref{stab-ctlr-cstrt-1}) can be rewritten in the following form of a network of Lur'e systems,
\begin{equation}
	\label{edge-rewrite}
	\begin{aligned}
		\dot{\tilde{z}}_2 &= \tA \tilde{z}_2 + \tB v, \\
		z &= ( U_{2} \otimes I_{n})\tilde{z}_2, \\
		v &= \Phi(z),
	\end{aligned}
\end{equation}
where $\tA=I_{N-1}\otimes A$, $\tB=\left(\Gamma U_{2}^T\right) \otimes (BK)$. 

The following theorem presents a sufficient condition for achieving the robust stabilization of (\ref{edge-dyn-3}) and equivalently the robust consensus of the initial MAS (\ref{mas}), and then how to design the consensus controller gain $K$.  

\begin{thm}
\label{rel-state-consensus-thm}
When $\Sigma_{1}$ and $\Sigma_{2}$ are not multipliers of identity matrices, the MAS (\ref{edge-dyn-3}) is robustly stabilized by the distributed stabilizing control law (\ref{stab-ctlr-cstrt-1}) and equivalently the robust consensus under relative-state constraints or uncertainties is achieved for the initial MAS (\ref{mas}) by the distributed controller (\ref{rel-state-consensus-ctlr-cstrt})
if there exist matrices $X \in \mathbb{R}^{n\times n}$, $Y \in \mathbb{R}^{m\times n}$ and $Z \in \mathbb{R}^{m\times m}$ such that the following LMI problem is feasible with $\epsilon > 0$,
\begin{equation}
	\label{rel-state-LMI}
	\begin{aligned}
		& \begin{bmatrix} \mathrm{sym}(AX+\lambda_{2}BY\Sigma_{2})+\epsilon X & \lambda_{2}BY+(\Sigma_{1}-\Sigma_{2})Z \\ \left(\lambda_{2}BY+(\Sigma_{1}-\Sigma_{2})Z\right)^T & -2Z \end{bmatrix} \preceq 0, \\
		& \begin{bmatrix} \mathrm{sym}(AX+\lambda_{N}BY\Sigma_{2})+\epsilon X & \lambda_{N}BY+(\Sigma_{1}-\Sigma_{2})Z \\ \left(\lambda_{N}BY+(\Sigma_{1}-\Sigma_{2})Z\right)^T & -2Z \end{bmatrix} \preceq 0, \\
		& X \succ 0, X ~\textrm{is diagonal}, \\
		& \begin{bmatrix} Z & X \\ X & \Psi^{-1} \end{bmatrix} \succeq 0, \\
		& \Psi \succ 0, \Psi \;\textrm{is diagonal}.		
	\end{aligned}
\end{equation}
Moreover, the controller gain $K$ is calculated by $K=YX^{-1}$. 
\end{thm}

\begin{IEEEproof}
Consider a Lyapunov function $V(\tilde{z}_2)=\tilde{z}_2^T\tP \tilde{z}_2$ where $\tP \triangleq I_{N-1}\otimes P$, $P \in \mathbb{R}^{n}$, $P \succ 0$. 
Taking the derivative of $V(\tilde{z}_2)$ gives us
\begin{equation*}
	\dot{V}(\tilde{z}_2) = \tilde{z}_2^T\left( \tP\tA+\tA^T\tP \right)\tilde{z}_2+2\tilde{z}_2^T\tP\tB v. 
\end{equation*}
Hence, for all $\epsilon > 0$ we have
\begin{equation*}
		\dot{V}(\tilde{z}_2)+\epsilon V(\tilde{z}_2) = \tilde{z}_2^T\left( \tP\tA+\tA^T\tP+\epsilon\tP \right)\tilde{z}_2+2\tilde{z}_2^T\tP\tB v.
\end{equation*} 
We now seek $\tP$ such that $\dot{V}(\tilde{z}_2)+\epsilon V(\tilde{z}_2) \leq 0$ as long as (\ref{rel-state-cstrt}) holds. 
Using the S-procedure \cite{Boyd:2004}, such $\tP$ exists if there exist $\psi_{1,1},\ldots,\psi_{1,n},\ldots,\psi_{M,1},\ldots,\psi_{M,n}$ which are non-negative such that 
\begin{align}
	\label{Vdot-1}
	& \dot{V}(\tilde{z}_2)+\epsilon V(\tilde{z}_2) \nonumber \\
	&  -\sum_{j=1}^{M}\sum_{k=1}^{n}{\psi_{j,k}(v_{j,k}-\sigma_{k,1}z_{j,k})(v_{j,k}-\sigma_{k,2}z_{j,k})} \leq 0. 
\end{align}
Let $\psi_{j,k}=\psi_{k} > 0 \; \forall \; j=1,\ldots,M$ and $\Psi=\mathrm{diag}\{\psi_{k}\}_{k=1,\ldots,n}$, then (\ref{Vdot-1}) is satisfied if 
\begin{align}
	\label{Vdot-2}
	& \dot{V}(\tilde{z}_2)+\epsilon V(\tilde{z}_2) -\sum_{j=1}^{M}{(v_{j}-\Sigma_{1}z_{j})^T\Psi(v_{j}-\Sigma_{2}z_{j})} \leq 0, \nonumber \\
	& \Leftrightarrow \begin{bmatrix} \tilde{z}_2 \\ v \end{bmatrix}^T \begin{bmatrix} \mathbb{P}_{1} & \mathbb{P}_{2} \\ \mathbb{P}_{2}^T & \mathbb{P}_{3} \end{bmatrix} 
	\begin{bmatrix} \tilde{z}_2 \\ v \end{bmatrix} \preceq 0 
	\Leftrightarrow \begin{bmatrix} \mathbb{P}_{1} & \mathbb{P}_{2} \\ \mathbb{P}_{2}^T & \mathbb{P}_{3} \end{bmatrix} \preceq 0,
\end{align}
where $\mathbb{P}_{1}=\tP\tA+\tA^T\tP+\epsilon\tP-I_{N-1}\otimes(\Psi \Sigma_{1}\Sigma_{2})$, $\mathbb{P}_{2}=\tP\tB+\frac{1}{2}U_{2}^T\otimes(\Psi(\Sigma_{1}+\Sigma_{2}))$, 
$\mathbb{P}_{3}=-I_{M}\otimes \Psi$.   

Subsequently, employing Schur complement \cite{Boyd:2004} to (\ref{Vdot-2}) results in 
$\mathbb{P}_{1}-\mathbb{P}_{2}\mathbb{P}_{3}^{-1}\mathbb{P}_{2}^T \preceq 0$, which is equivalent to
\begin{align}
	\label{Vdot-3}
	& I_{N-1}\otimes(A^TP+PA+\epsilon P-\Psi \Sigma_{1}\Sigma_{2})  \nonumber \\
	& +\frac{1}{2}\Gamma\otimes\mathrm{sym}(PBK(\Sigma_{1}+\Sigma_{2})) 
	 +\Gamma^2\otimes(PBK\Psi^{-1}K^TB^TP) \nonumber \\ 
	& +\frac{1}{4}I_{N-1}\otimes\left[\Psi(\Sigma_{1}+\Sigma_{2})^2\right] \preceq 0.
\end{align}
Since $\Gamma$ is diagonal, (\ref{Vdot-3}) is equivalent to 
\begin{align*}
	& A^TP+PA+\epsilon P-\Psi \Sigma_{1}\Sigma_{2}+\lambda_{k}^2PBK\Psi^{-1}K^TB^TP \nonumber \\
	& +\frac{1}{2}\lambda_{k}\mathrm{sym}(PBK(\Sigma_{1}+\Sigma_{2}))+\frac{1}{4}\Psi(\Sigma_{1}+\Sigma_{2})^2  \preceq 0,  \nonumber \\
	\Leftrightarrow & A^TP+PA+\epsilon P+\lambda_{k}^2PBK\Psi^{-1}K^TB^TP \nonumber \\
	& +\frac{1}{2}\lambda_{k}\mathrm{sym}(PBK(\Sigma_{1}+\Sigma_{2}))+\frac{1}{4}\Psi(\Sigma_{1}-\Sigma_{2})^2  \preceq 0.
\end{align*}
Next, denote $X \triangleq P^{-1}$ and multiply $X$ both to the left and to the right of the equation above, we obtain
\begin{align}
	\label{Vdot-5}
	& XA^T+AX+\epsilon X+\lambda_{k}^2BK\Psi^{-1}K^TB^T \nonumber \\
	& +\frac{1}{2}\lambda_{k}\mathrm{sym}(BK(\Sigma_{1}+\Sigma_{2})X)+\frac{1}{4}X\Psi(\Sigma_{1}-\Sigma_{2})^2X \preceq 0.
\end{align} 
If $X$ is diagonal then (\ref{Vdot-5}) is equivalent to
\begin{align}
	\label{Vdot-6}
	& \mathrm{sym}(AX+\lambda_{k}BY\Sigma_{2})+\epsilon X+\frac{1}{2}[\lambda_{k}BY+(\Sigma_{1}-\Sigma_{2})Z] \nonumber \\
	& \times Z^{-1}[\lambda_{k}BY+(\Sigma_{1}-\Sigma_{2})Z]^T \preceq 0,	
\end{align} 
where $Y \triangleq KX$, $Z=\frac{1}{2}X^2\Psi$.  
Then using Schur complement again with (\ref{Vdot-6}) leads to 
\begin{equation}
	\begin{bmatrix} \mathrm{sym}(AX+\lambda_{k}BY\Sigma_{2})+\epsilon X & \lambda_{k}BY+(\Sigma_{1}-\Sigma_{2})Z \\ \left(\lambda_{k}BY+(\Sigma_{1}-\Sigma_{2})Z\right)^T & -2Z \end{bmatrix} \preceq 0,		
\end{equation}
for all $k=2,\ldots,N$. 
Since $\lambda_{2} \leq \lambda_{3},\ldots,\lambda_{N-1} \leq \lambda_{N}$, we can represent $\lambda_{i},i=3,\ldots,N-1$ as convex combinations of $\lambda_{2}$ and $\lambda_{N}$. 
Thus, we derive (\ref{rel-state-LMI}). 
\end{IEEEproof}

\begin{thm}
\label{rel-state-consensus-thm-1}
Suppose that $\Sigma_{1} =\sigma_{1}I_{n}$ and $\Sigma_{2} =\sigma_{2}I_{n}$ then the MAS (\ref{edge-dyn-3}) is robustly stabilized by the distributed stabilizing control law (\ref{stab-ctlr-cstrt-1}) and equivalently the robust consensus under relative-state constraints or uncertainties is achieved for the initial MAS (\ref{mas}) by the distributed controller (\ref{rel-state-consensus-ctlr-cstrt})
if there exist matrices $X \in \mathbb{R}^{n\times n}$, $Y \in \mathbb{R}^{m\times n}$ and $Z \in \mathbb{R}^{m\times m}$ such that the following LMI problem is feasible with $\epsilon > 0$,
\begin{equation}
	\label{rel-state-LMI-1}
	\begin{aligned}
		& \begin{bmatrix} \mathrm{sym}(AX+\sigma_{2}\lambda_{2}BY)+\epsilon X & \lambda_{2}BY+(\sigma_{1}-\sigma_{2})Z \\ \left(\lambda_{2}BY+(\sigma_{1}-\sigma_{2})Z\right)^T & -2Z \end{bmatrix} \preceq 0, \\
		& \begin{bmatrix} \mathrm{sym}(AX+\sigma_{2}\lambda_{N}BY)+\epsilon X & \lambda_{N}BY+(\sigma_{1}-\sigma_{2})Z \\ \left(\lambda_{N}BY+(\sigma_{1}-\sigma_{2})Z\right)^T & -2Z \end{bmatrix} \preceq 0, \\
		& X \succ 0,  \\
		& \begin{bmatrix} Z & X \\ X & \Psi^{-1} \end{bmatrix} \succeq 0, \\
		& \Psi \succ 0, \Psi \;\textrm{is diagonal}.		
	\end{aligned}
\end{equation}
Moreover, the controller gain $K$ is calculated by $K=YX^{-1}$. 
\end{thm}

\begin{IEEEproof}
Consider the same Lyapunov function as in the proof of Theorem \ref{rel-state-consensus-thm}. 
Then all steps until (\ref{Vdot-5}) are also true in this scenario. 
Accordingly, substituting $\Sigma_{1} =\sigma_{1}I_{n}$ and $\Sigma_{2} =\sigma_{2}I_{n}$ into (\ref{Vdot-5}) gives us
\begin{align}
	\label{Vdot-7}
	& \mathrm{sym}(AX+\sigma_{2}\lambda_{k}BY)+\epsilon X+\frac{1}{2}[\lambda_{k}BY+(\sigma_{1}-\sigma_{2})Z] \nonumber \\
	& \times Z^{-1}[\lambda_{k}BY+(\sigma_{1}-\sigma_{2})Z]^T \preceq 0,	
\end{align} 
where $Y \triangleq KX$, $Z=\frac{1}{2}X\Psi X$.  
Then using Schur complement again with (\ref{Vdot-7}) and notes that $\lambda_{i},i=3,\ldots,N-1$ can be represented as convex combinations of $\lambda_{2}$ and $\lambda_{N}$, 
we obtain (\ref{rel-state-LMI-1}). 
\end{IEEEproof}

\begin{rem}
Recently, there are several existing researches, e.g. \cite{Franceschelli:2013}, \cite{Tran:2015}, which propose different distributed methods to approximate the whole eigen-spectrum of the Laplacian matrix. These methods can be employed to estimate $\lambda_{2}$ and $\lambda_{N}$ before solving the LMI problems (\ref{rel-state-LMI}), (\ref{rel-state-LMI-1}).  
As a result, we can solve (\ref{rel-state-LMI}) and (\ref{rel-state-LMI-1}) in a distributed fashion.   
\end{rem}

\begin{rem}
The difference between Theorem \ref{rel-state-consensus-thm-1} and Theorem \ref{rel-state-consensus-thm} is that the variable $X$ in (\ref{rel-state-LMI-1}) is not required to be diagonal 
while that in (\ref{rel-state-LMI}) is. Therefore, if $\Sigma_{1}$ and $\Sigma_{2}$ are multipliers of identity matrices, i.e., the upper and lower sector slopes for relative state constraints or uncertainties of all agents are the same then the associated LMI problem is less conservative and hence its feasibility is increased.   
\end{rem}

\begin{rem}
As stated in the introduction, our method to derive LMI problems (\ref{rel-state-LMI}) and (\ref{rel-state-LMI-1}) for the Lur'e network (\ref{edge-rewrite}) 
is more general than the method for a similar Lur'e network in \cite{Takaba:2014}.  
On the other hand, the problem setting in this paper is in a different form of Lur'e networks with that in \cite{Zhang:2014j}. Therefore, the obtained results are not similar.  
More specifically, \cite{Zhang:2014j} uses a linear cooperative input and another nonlinear input with a different input matrix $E$ satisfying the incrementally passive or incrementally sector-bounded condition, which is less general than our sector-bounded condition (\ref{rel-state-cstrt}). 
\end{rem}

\section{Numerical Examples}
\label{app}

\subsection{Practical Consensus of Mobile Robots}

\begin{figure}[htpb!]
		\centering
		\includegraphics[scale=0.3]{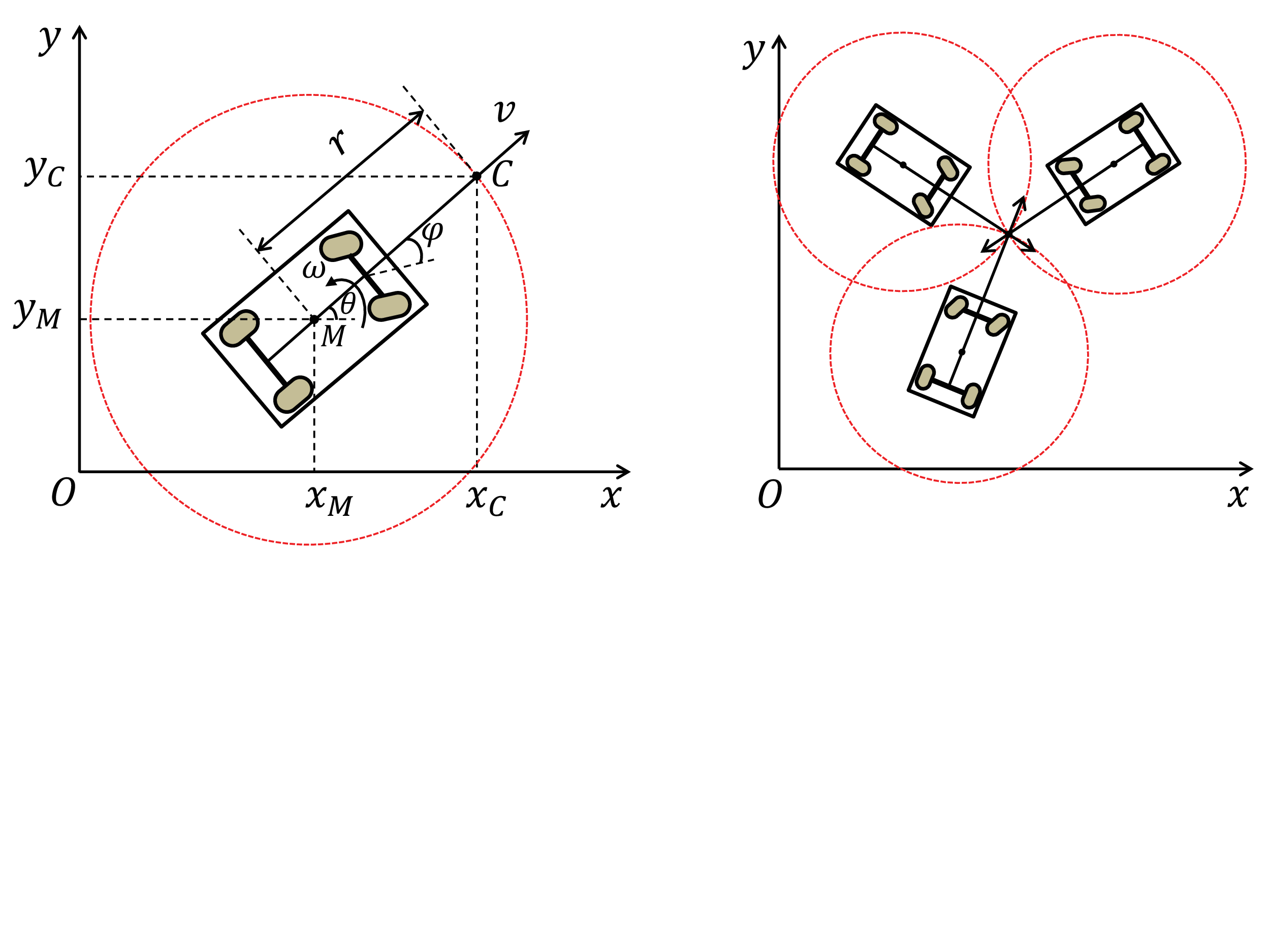}
		\vspace{-25mm}
		\caption{Demonstration of variables related to a $4$-wheel robot and the concept of practical robot consensus.}
		\label{robot_concept}
	\end{figure}

Consider a group of $N$ identical $4$-wheel robots with front-wheel steering. 
The variables and parameters of each robot are illustrated in Figure \ref{robot_concept} where the center of mass is denoted by $M_{i}$ whose position in a given coordinate $(O,x,y)$ is represented by 
$(x_{Mi},y_{Mi})$. The rotation and steering angles are denoted by $\theta_{i}$ and $\varphi_{i}$, respectively. 
Accordingly, $\omega_{i}$ and $v_{i}$ represent the angular and longitudinal velocity. 
To take into account practical factors such as robots' dimensions and collision avoidance, 
we shall investigate the consensus of the robots' heading points $C_{i}$ instead of $M_{i}$. 
This practical consensus concept is demonstrated in Figure \ref{robot_concept}.

Next, the robot's model in term of the coordinates of $C_{i}$ is as follows, 
\begin{align}
	\label{2w-steering-robot-ss}
	\dot{x}_{i} = M_{i}[v_{i},\omega_{i}]^T, i=1,\ldots,N,
\end{align}
where $x_{i}=[x_{Ci},y_{Ci}]^T; ~ M_{i}=\begin{bmatrix} \cos\theta_{i} & -r\sin\theta_{i} \\ \sin\theta_{i} & r\cos\theta_{i} \end{bmatrix}.$ 
Denote $u_{i} \triangleq M_{i}[v_{i},\omega_{i}]^T$ then each robot can be represented by a set of two integrators 
$\dot{x}_{i} = u_{i}, i=1,\ldots,N.$  
We have
\begin{equation*}
	\begin{bmatrix} v_{i} \\ \omega_{i} \end{bmatrix} = M_{i}^{-1}u_{i} = \begin{bmatrix} \cos\theta_{i} & \sin\theta_{i} \\ -\sin\theta_{i}/r & \cos\theta_{i}/r \end{bmatrix} u_{i} 
	\triangleq \begin{bmatrix} \tilde{u}_{i,1} \\ \tilde{u}_{i,2} \end{bmatrix}. 
\end{equation*}
Therefore, the real control inputs $v_{i}$ and $\varphi_{i}$ to each robot are computed by
\begin{equation}
	\label{real-input}
	v_{i}=\tilde{u}_{i,1}; ~ \varphi_{i}=\arctan \frac{\omega_{i}}{v_{i}}=\arctan \frac{\tilde{u}_{i,2}}{\tilde{u}_{i,1}}.
\end{equation}
Consequently, we consider the constraint  $\|x_{i}-x_{j}\|_{\infty} \leq \alpha$ on relative states of connected robots,  
which implies that the communication range between robots is limited to $\sqrt{2}\alpha$. 
This is indeed a robust consensus problem under relative-state constraints within our framework. 

\begin{figure}[htpb!]
		\centering
		\includegraphics[scale=0.35]{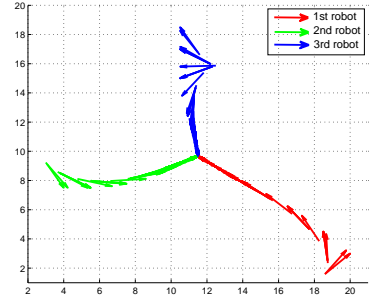}
		\caption{Trajectories of $4$-wheel robots reaching consensus.}
		\label{robot_consensus}
	\end{figure}

\begin{figure}[htpb!]
		\centering
		\includegraphics[scale=0.35]{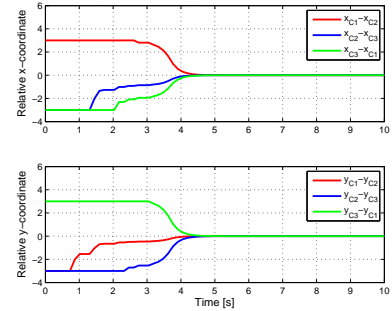}
		\caption{Bounded relative coordinates of heading points of $4$-wheel robots.}
		\label{robot_xy_relstate}
	\end{figure}

Employing Theorem \ref{rel-state-consensus-thm-1}, we solve the LMI problem (\ref{rel-state-LMI-1}) with $A=0,B=1$ and obtain $\frac{-\epsilon}{\lambda_{2}\sigma_{2}}<K<0.$ 
In the simulation, we set $r=2$ [dm], $\alpha=3$ [dm], $\epsilon=0.4$, and $\tG$ is a full graph, then choose $K=-0.1$ since $\lambda_{2}=3$. 
The simulation result in Figure \ref{robot_consensus} then confirms that the consensus among robots is achieved even though there is a constraint on relative state exchange of robots, 
where the arrows represent the vectors $\overrightarrow{M_{i}C_{i}}$ of robots. 
Moreover, Figure \ref{robot_xy_relstate} shows that the exchanged relative states of robots always satisfy the given bounded constraint.

\subsection{Consensus of Oscillator Networks}

\begin{figure}[ht!]
		\centering
		\includegraphics[scale=0.35]{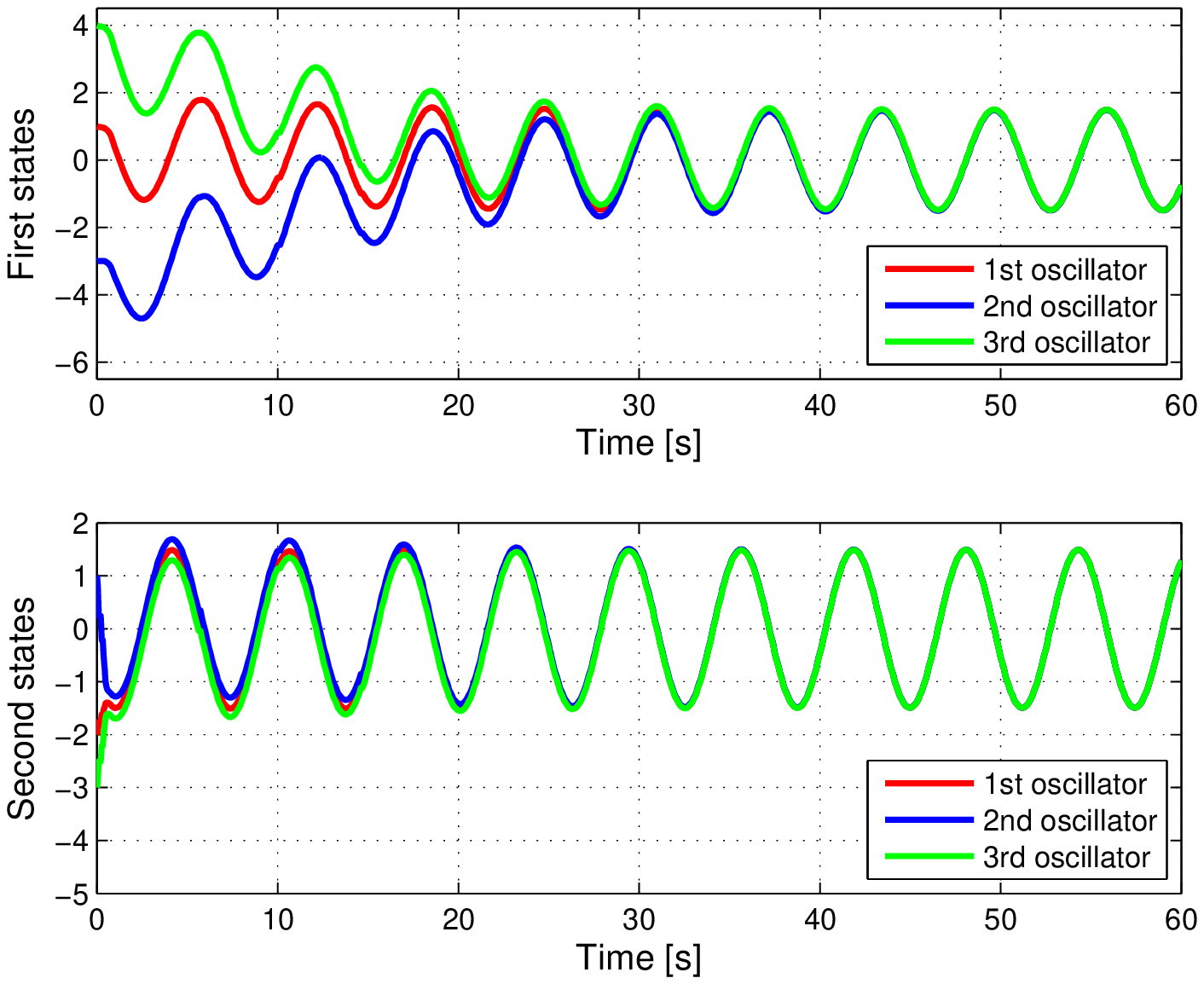}
		\caption{Synchronization of oscillators under relative state constraints.}
		\label{osc-rel-state-cst-1}
		\includegraphics[scale=0.35]{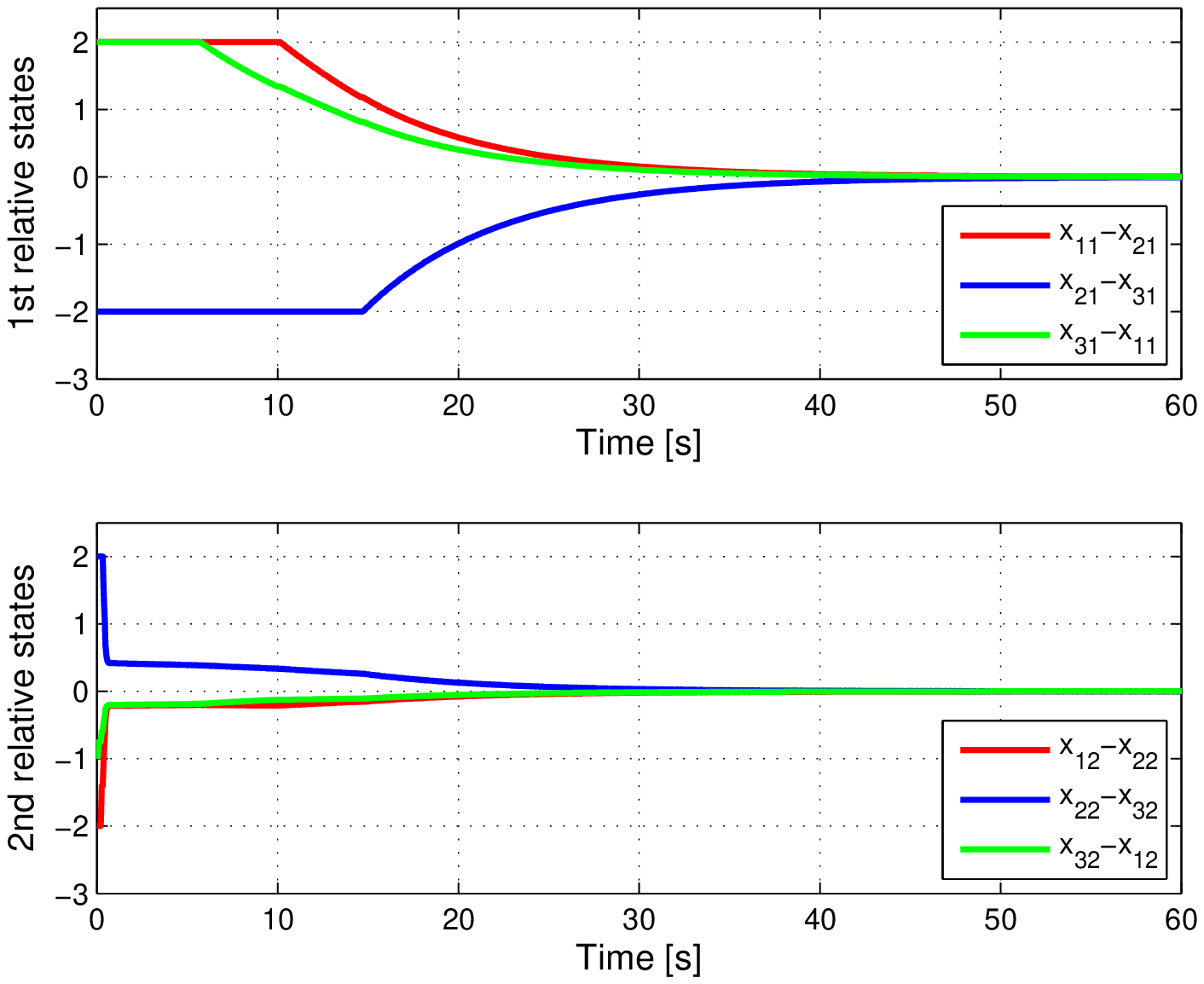}
		\caption{Bounded relative states of oscillators.}
		\label{osc-rel-state-cst-2}		
	\end{figure}

To further illustrate the proposed approach, we consider a consensus problem in a network of $3$ identical linear oscillators with the following model,
\begin{equation}
	\label{lin-osc}
	\dot{x}_k=Ax_k+Bu_k,k=1,2,3,
\end{equation}
where
\begin{equation}
	A=\begin{bmatrix} 0 & 1 \\ -1 & 0 \end{bmatrix}, B=\begin{bmatrix} 0 \\ 1 \end{bmatrix},
\end{equation}
and the initial conditions of the oscillators are $(1,-2)$; $(-3,1)$; $(4,-3)$, respectively. 
We then assume that $\tG$ is a full graph and the exchanged relative states among agents are bounded in $[-2,2]$. 
With $\epsilon=0.1$, solving the LMI problem (\ref{rel-state-LMI-1})  gives us $K=[-2.3825,-20.6800]$. 
Consequently, Figure \ref{osc-rel-state-cst-1} reveals that the oscillators exhibit synchronized oscillations whereas Figure \ref{osc-rel-state-cst-2} shows that 
the relative states of oscillators satisfy the bounded constraints, i.e., the robust consensus is achieved.

\begin{figure}[ht!]
		\centering
		\includegraphics[scale=0.35]{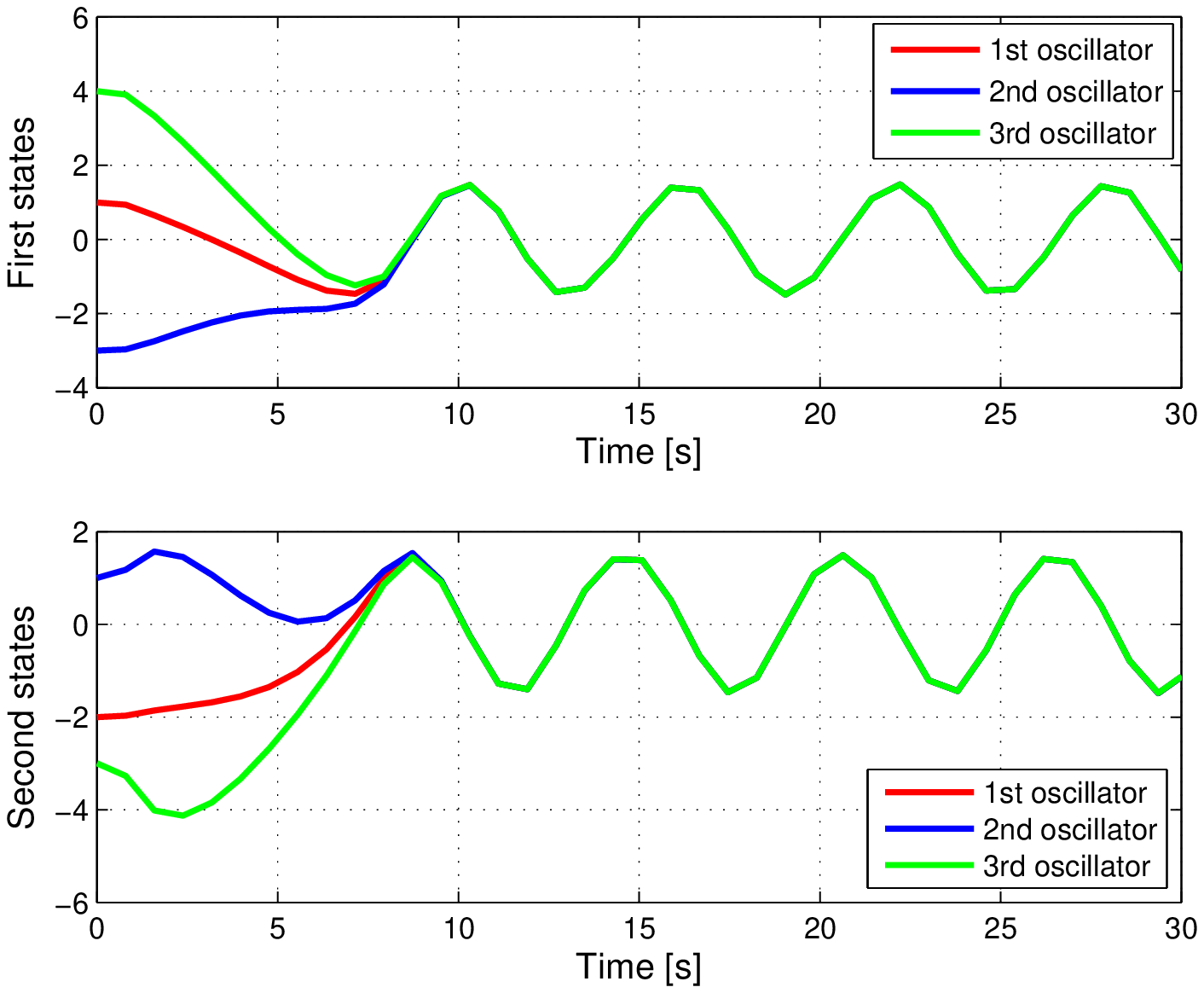}
		\caption{Synchronization of oscillators under relative state uncertainties.}
		\label{osc-rel-state-unc-1}
		\includegraphics[scale=0.35]{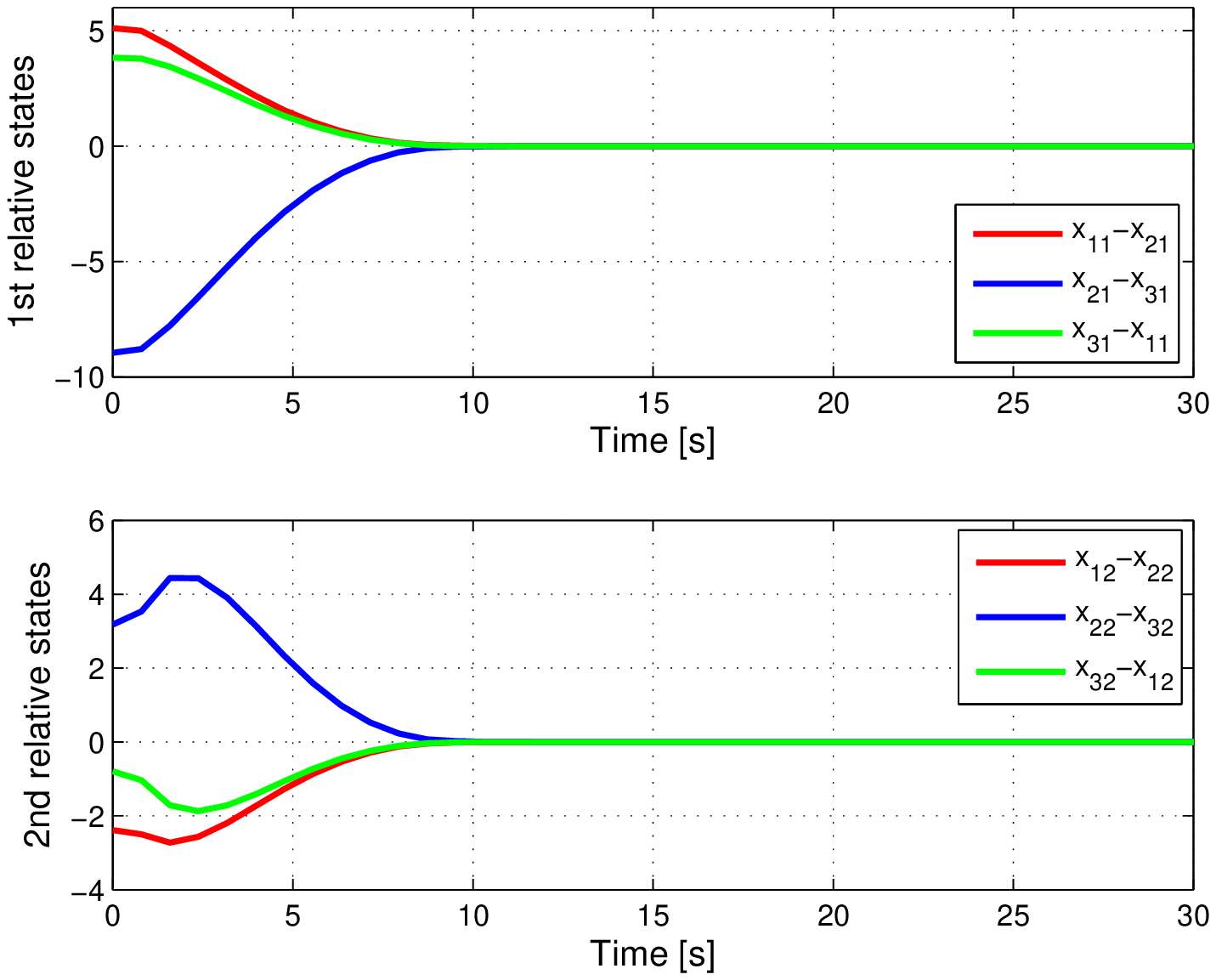}
		\caption{Relative states of oscillators.}
		\label{osc-rel-state-unc-2}		
	\end{figure}

Next, we demonstrate the effectiveness of our approach in the scenario that the exchanged relative states among agents contain some uncertainties that results in $\Sigma_{1}=0.7I_{2}$ 
and $\Sigma_{2}=1.3I_{2}$. Then we solve the LMI problem (\ref{rel-state-LMI-1}) to obtain $K=[-1.1309,-2.2191]$. In the simulation, we randomly generate those uncertainties in the interval $[0.7,1.3]$. 
We then observe that the synchronization of oscillators are achieved for any uncertainties in the given range. 
Particularly, Figure \ref{osc-rel-state-unc-1}--\ref{osc-rel-state-unc-2} display the oscillator network's responses 
for a specific case where the uncertainties on two relative states of oscillators are $1.2789$ and $0.7946$.

\section{Conclusion}
\label{sum}

An approach has been proposed in this paper to analyze and synthesize distributed global robust consensus controllers for general linear leaderless MASs under relative state constraints or uncertainties 
with the following appealing features. 
First, it is available for a broader class of MASs and for constraints or uncertainties described by a sector-bounded condition 
which is more general than that in the existing researches. 
Second, it shows that the global robust consensus design with relative state constraints or uncertainties is equivalent to a robust stabilizing design with state constraints or uncertainties 
of a transformed MAS. 
Third, a sufficient condition for global robust consensus and the global robust consensus controller gain are derived from the solutions of a distributed convex LMI problem 
which is less conservative than in other studies.


%


\section*{Acknowledgment}


The author would like to thank Toyota Technological Institute for its supports.

\ifCLASSOPTIONcaptionsoff
  \newpage
\fi



%

\bibliographystyle{IEEEtran}
\bibliography{IEEEabrv,References}

%
%

%

%
%
%
%
%
%
%
%
%
\end{document}